\documentstyle{article}

\def\mj{{\mathbf{1}}}
\def\pl{\!+\!}
\def\mn{\!-\!}
\def\cirk{\,{\raisebox{.3ex}{\tiny $\circ$}}\,}
\def\cirkk{\,{\raisebox{.1ex}{\tiny $\circ$}}\,}
\def\prop#1#2{\vspace{2ex} \noindent{\sc #1.} {\it #2} \par \vspace{2ex}}
\def\dkz{\noindent{\sc Proof. }}
\def\qed{\hfill $\dashv$}

\begin{document}

\title{Isomorphic Formulae\\ in Classical Propositional Logic}
\author{\small {\sc Kosta Do\v sen} and {\sc Zoran Petri\' c}
\\[1ex]
{\small Mathematical Institute, SANU}\\[-.5ex]
{\small Knez Mihailova 36, p.f.\ 367, 11001 Belgrade,
Serbia}\\[-.5ex]
{\small email: \{kosta, zpetric\}@mi.sanu.ac.rs}}
\date{}
\maketitle

\begin{abstract}
\noindent Isomorphism between formulae is defined with respect to
categories formalizing equality of deductions in classical
propositional logic and in the multiplicative fragment of
classical linear propositional logic caught by proof nets. This
equality is motivated by generality of deductions.
Characterizations are given for pairs of isomorphic formulae,
which lead to decision procedures for this isomorphism.
\end{abstract}


\vspace{.5ex}


\section{Introduction}
Isomorphism between formulae should be an equivalence relation
stronger than mutual implication. This is presumably the relation
underlying the relation that holds between propositions that have
the same meaning just because of their logical form. Any
propositions that are instances, with the same substitution, of
isomorphic formulae would have the same meaning, which presumably
need not be the case for formulae that are just equivalent, i.e.\
which just imply each other.

One may try to characterize isomorphic formulae by looking only
into the inner structure of formulae . This is the way envisaged
by Carnap and Church (see \cite{C47}, Sections 14-15, where
related work by Quine and C.I.\ Lewis is mentioned, \cite{Ch50},
\cite{A98}, Section~2, and references therein).

Another way is to try to characterize isomorphism between formulae
by looking also at the outer structure in which formulae are to be
found. This outer structure may be a deductive structure, which is
characterized in terms of categories in categorial proof theory.
The categories we need are syntactical: their objects are formulae
and their arrows are deductions.

Isomorphism between formulae may then be understood exactly as
isomorphism between objects is understood in category theory. The
formulae $A$ and $B$ are isomorphic when there is a deduction $f$,
i.e.\ arrow, from $A$ to $B$, and another deduction $g$ from $B$
to $A$, such that $f$ composed with $g$ is equal to the identity
deduction from $A$ to $A$, while $g$ composed with $f$ is equal to
the identity deduction from $B$ to $B$. This analysis of
isomorphism presupposes a notion of equality between deductions,
which is formalized in our syntactical categories. (\emph{Equality
between deductions} stands here for what we and other authors have
called elsewhere \emph{identity of proofs}; see \cite{D03},
\cite{DP04}, Sections 1.3-1.4, and references therein.)
Characterizing this notion of equality between deductions is the
main task of categorial proof theory, and of general proof theory.

That $A$ and $B$ are isomorphic means here intuitively that they
function in the same manner in deductions. In a deduction one can
replace one by the other, either as premise or as conclusion, so
that nothing is lost, nor gained. The replacements, which are made
by composing our deduction with the deductions $f$ and $g$, are
such that they enable us to return to our original deduction by
further composing with $g$ and $f$, since $f$ composed with $g$
and $g$ composed with $f$ are identity deductions, and hence may
be cancelled. (For a view concerning isomorphic formulae like the
one presented here, and its relationship with propositional
identity, see \cite{D95}, Section~9, and \cite{D99}, Section~5.)

The study of isomorphic formulae first started in intuitionistic
logic, for which it is widely believed that we have a solid
nontrivial notion of equality of deductions. This notion is
characterized either in terms of the typed lambda calculus (via
the Curry-Howard correspondence), or in terms of categories based
on cartesian closed categories (these characterizations may be
equivalent; see \cite{LS86}).

A result exists in this area for the
implication-conjunction-$\top$ fragment of intuitionistic logic
(see \cite{S81}). As far as we know, the latest advances
concerning the still open problem of characterizing formulae
isomorphic in the whole of intuitionistic propositional logic
(which is related to Tarski's high-school algebra problem; see
\cite{BL93}) were made in \cite{B02} and \cite{FDB06}. There is a
related result characterizing isomorphic formulae in the analogous
multiplicative fragment of linear logic, which corresponds to
symmetric monoidal closed categories, and is common to classical
and intuitionistic linear logic (see \cite{S93} and \cite{DP97}).

The problem of characterizing isomorphic formulae was not
approached up to now in classical propositional logic, and the
results we are going to present here cover this logic. They cover
also a fragment of classical linear propositional logic. To be
able to approach our problem, and to obtain significant results,
we need for the logics we want to cover a plausible and nontrivial
notion of equality of arrows in categories formalizing equality of
deductions in these logics. A consensus for classical linear
propositional logic may be found around the multiplicative
fragment of that logic caught by proof nets, which leads to
notions of category closely related to star-autonomous category
(see \cite{S89}, \cite{DP07} and references therein).

For classical propositional logic, it is on the contrary widely
believed that no nontrivial notion of category would do the job.
It is believed that no nontrivial notion of Boolean category may
be found. This is indeed the case if one wants these Boolean
categories to be cartesian closed (see \cite{D03}, Section~5,
\cite{DP04}, Section 14.3, and references therein). But, whereas
on the level of theorems classical logic is an extension of
intuitionistic logic, it is not clear that this should be so at
the level of deductions and of their equality.

If one does not require that Boolean categories be cartesian
closed categories, and bases equality of deductions in Boolean
categories on coherence results analogous to those available for
classical linear propositional logic, a nontrivial notion may
arise. The coherence results in question are categorial results
analogous to the classical coherence result of Kelly and Mac Lane
for symmetric monoidal closed categories (see \cite{KML71}). They
reduce equality of arrows in the syntactical category to equality
of arrows in a graphical model category.

Such a nontrivial notion of Boolean category may be found in
\cite{DP04} (Chapter 14), and Section~4 of the present paper deals
with isomorphism of formulae engendered by that particular notion.
Section~3 of the paper deals with isomorphism of formulae in
classical and classical linear propositional logic different from
the notion of Section~4. That other notion, which involves
graphical model categories quite like those of Kelly and Mac Lane,
is motivated by generality of deductions (as suggested by
\cite{L69}). The notion of Section~4 may also be understood as
involving generality up to a point, but the notion of Section~3
does so more consistently. Both notions are however analogous in
that they base equality of deductions on equality of arrows in
some graphical model category.

The main results of the paper in Sections 3 and 4 are based on
some preliminary elementary results concerning classical
propositional logic, which are established in Section~2, and
occasionally later in the paper. Although these results are not
difficult to reach, when they are combined with more advanced
results, such as those that may be found in \cite{DP04}, they give
a complete characterization of isomorphic formulae in classical
and classical linear propositional logic. These characterizations
are such that they easily lead to decision procedures for the
isomorphisms in question.

This paper is devoted to the problem of characterizing pairs of
isomorphic formulae. A related, but different, problem involving
isomorphism is to characterize arrows that are isomorphisms, in
categories formalizing equality of deductions. (This problem for
the conjunction-$\top$ fragment of classical or intuitionistic
logic is dealt with in \cite{DP99}.) Although we have not dealt
with that second problem explicitly in this paper, a solution for
it may easily be inferred from our results.

We will however not dwell on that, in order not to overburden the
text with categorial matters. In the whole paper we try to keep
the presence of categories to a minimum, and give more prominence
to elementary, easily understandable, logical facts. To appreciate
the full import of our results the reader should however be
acquainted up to a point with certain notions covered in
particular by \cite{DP04} (and which we cannot possibly expose all
here).

\section{$\wedge\vee$ and $\neg\!\wedge\!\vee$-equivalences}
Let $\cal L$ be a propositional language generated out of an
infinite set of propositional letters, which we call simply
\emph{letters}, with the nullary connectives $\top$ and $\bot$,
the unary connective $\neg$, and the binary connectives $\wedge$
and $\vee$. We use $p,q,r,\ldots,$ sometimes with indices, for
letters, and $A,B,C,\ldots,$ sometimes with indices, for the
formulae of $\cal L$. Let ${\cal L}_{\wedge,\vee}$, ${\cal
L}_{\neg,\wedge,\vee}$ and ${\cal L}_{\top,\bot,\wedge,\vee}$ be
the propositional languages defined as $\cal L$ but with only the
connectives in the subscripts.

We envisage also propositional languages extending $\cal L$, in
which we have moreover the binary connectives of equivalence
$\leftrightarrow$ and implication $\rightarrow$. We may imagine
that these two additional connectives are defined in $\cal L$, and
we shall not introduce special names for these extended languages.

Let a $\wedge\vee$\emph{-equivalence} be a formula
${A\leftrightarrow B}$ where $A$ and $B$ are formulae of ${\cal
L}_{\wedge,\vee}$. Consider the formal system ${\cal
S}_{\wedge,\vee}$ whose axioms are the $\wedge\vee$-equivalences
of the following forms:
\begin{tabbing}
\centerline{${A\leftrightarrow A}$,}
\\[1ex]
\centerline{$((A\wedge B)\wedge C)\leftrightarrow(A\wedge (B\wedge
C))$,\hspace{3em}$((A\vee B)\vee C)\leftrightarrow(A\vee (B\vee
C))$,}
\\[1ex]
\centerline{$(A\wedge B)\leftrightarrow(B\wedge
A)$,\hspace{8.4em}$(A\vee B)\leftrightarrow(B\vee A)$,}
\end{tabbing}
and whose theorems are the $\wedge\vee$-equivalences obtained
starting from these axioms with the following rules:
\[
\hspace{2.5em}\frac{A\leftrightarrow B}{B\leftrightarrow
A}\hspace{11em}\frac{A\leftrightarrow B \quad\quad
B\leftrightarrow C}{A\leftrightarrow C}
\]
\[
\frac{A\leftrightarrow B \quad\quad C\leftrightarrow D}{(A\wedge
C)\leftrightarrow (B\wedge D)}\hspace{8.4em}\frac{A\leftrightarrow
B \quad\quad C\leftrightarrow D}{(A\vee C)\leftrightarrow (B\vee
D)}
\]

A formula is \emph{diversified} when every letter occurs in it not
more than once. A $\wedge\vee$-equivalence is \emph{diversified}
when $A$ and $B$ are diversified.

Assume that ${\mbox{\it let}}A$ is the set of letters occurring in
$A$, and let $A^p_B$ be the result of substituting the formula $B$
for every occurrence of $p$ in $A$. We can prove the following
lemmata.

\prop{Lemma~1}{Assume that $A$ is a diversified formula of ${\cal
L}_{\wedge,\vee}$, that $B$ is a subformula of $A$, and that
${\mbox{\it let}}A\mn{\mbox{\it let}}B=\{p_1,\ldots,p_n\}$ (where
$\{p_1,\ldots,p_n\}$ is empty if $n=0$). Then there is a sequence
$S_1,\dots, S_n$, where $S_i\in\{\top,\bot\}$, such that
$A^{p_1\ldots p_n}_{S_1\ldots S_n}\leftrightarrow B$ is a
tautology.}

\dkz We proceed by induction on $n$. If $n=0$, then $B$ is $A$,
and $A\leftrightarrow B$ is of course a tautology.

If $n=k\pl 1$, and $A$ is of the form ${C\wedge D}$ with $B$ a
subformula of $C$, then for ${\mbox{\it let}}C\mn{\mbox{\it
let}}B=\{q_1,\ldots,q_m\}$, where
$\{q_1,\ldots,q_m\}\cup\{r_1,\ldots,r_l\}=\{p_1,\ldots,p_n\}$ for
$m\geq 0$ and $l\geq 1$, we have by the induction hypothesis that
for some $S_1,\dots, S_m$ the formula $C^{q_1\ldots
q_m}_{S_1\ldots S_m}\leftrightarrow B$ is a tautology. Hence
$A^{q_1\ldots q_m r_1\ldots r_l}_{S_1\ldots
S_m\top\ldots\top}\leftrightarrow B$ is a tautology too.

We proceed analogously when $B$ is a subformula of $D$, or when
$A$ is of the form ${C\vee D}$. In the latter case we substitute
$\bot,\ldots,\bot$ for $r_1,\ldots,r_l$.\qed

\prop{Lemma~2}{If $B$ is a diversified formula of ${\cal
L}_{\wedge,\vee}$ that has a letter $q$ not in the formula $A$ of
$\cal L$, then ${A\leftrightarrow B}$ is not a tautology.}

\dkz By Lemma~1, for some $S_1,\dots, S_n$, where
$S_i\in\{\top,\bot\}$, we have that $B^{p_1\ldots p_n}_{S_1\ldots
S_n}\leftrightarrow q$ is a tautology. On the other hand,
$A^{p_1\ldots p_n}_{S_1\ldots S_n}\leftrightarrow q$ cannot be a
tautology. So $A^{p_1\ldots p_n}_{S_1\ldots S_n}\leftrightarrow
B^{p_1\ldots p_n}_{S_1\ldots S_n}$ is not a tautology, and hence
${A\leftrightarrow B}$ is not a tautology.\qed

\vspace{2ex}

For $A$ a diversified formula of $\cal L_{\wedge,\vee}$ we say
that the letters $p$ and $q$ are \emph{conjunctively joined} in
$A$ when $A$ has a subformula ${P\wedge Q}$ or ${Q\wedge P}$ such
that $p$ is in $P$ and $q$ is in $Q$, and we say that $p$ and $q$
are \emph{directly} conjunctively joined in $A$ when $A$ has a
subformula ${P\wedge Q}$ or ${Q\wedge P}$ such that $p$ is in $P$
and $q$ is in $Q$, no subformula of $P$ containing $p$ is a
disjunction, and no subformula of $Q$ containing $q$ is a
disjunction. We define analogously disjunctively and directly
disjunctively joined formulae (we just replace $\wedge$ by $\vee$
and $\vee$ by $\wedge$). We can prove the following proposition.

\prop{Proposition~1}{A diversified $\wedge\vee$-equivalence is a
tautology iff it is a theorem of ${\cal S}_{\wedge,\vee}$.}

\dkz It is clear that every theorem of ${\cal S}_{\wedge,\vee}$ is
a tautology. (This is established by induction on the length of
proof in ${\cal S}_{\wedge,\vee}$.) For the converse, we suppose
that the diversified $\wedge\vee$-equivalence ${A\leftrightarrow
B}$ is a tautology, and we proceed by induction on the number $n$
of connectives in $A$. By Lemma~2, this number must also be the
number of connectives in $B$.

If $n=0$, then, by Lemma~2, the formulae $A$ and $B$ must both be
a letter $p$, and $p\leftrightarrow p$ is an axiom of ${\cal
S}_{\wedge,\vee}$. If $n=k\pl 1$, then $A$ has a subformula either
of the form ${p\wedge q}$ or of the form ${p\vee q}$. Suppose
${p\wedge q}$ is a subformula of $A$, and consider how $p$ and $q$
are joined in $B$.

(1) It is impossible that $p$ and $q$ be disjunctively joined in
$B$. Suppose they are. By Lemma~1, for some $S_1,\dots, S_k$,
where $S_i\in\{\top,\bot\}$, we have that $A^{r_1\ldots
r_k}_{S_1\ldots S_k}\leftrightarrow (p\wedge q)$ is a tautology.
On the other hand, by using Lemma~2 we infer that $B^{r_1\ldots
r_k}_{S_1\ldots S_k}\leftrightarrow C$ can be a tautology only if
either ${C\leftrightarrow (p\vee q)}$, or ${C\leftrightarrow p}$,
or ${C\leftrightarrow q}$, or ${C\leftrightarrow\top}$, or
${C\leftrightarrow\bot}$, is a tautology. So $A^{r_1\ldots
r_k}_{S_1\ldots S_k}\leftrightarrow B^{r_1\ldots r_k}_{S_1\ldots
S_k}$ is not a tautology, which contradicts the assumption that
${A\leftrightarrow B}$ is a tautology.

(2) It is also impossible that $p$ and $q$ be conjunctively joined
in $B$, but not directly. Otherwise, we would have in $B$ a
subformula ${P\wedge Q}$ or ${Q\wedge P}$ with $p$ in $P$ and $q$
in $Q$, and a subformula $C\vee D$ of $P$ with $p$ either in $C$
or in $D$; we need not consider separately the analogous case when
$C\vee D$ is a subformula of $Q$ with $q$ either in $C$ or in $D$.
Then ${A^p_\bot\leftrightarrow A'}$ and ${B^p_\bot\leftrightarrow
B'}$ are tautologies for $A'$ a formula of $\cal L$ without $q$
and $b'$ a diversified formula of ${\cal L}_{\wedge,\vee}$ with
$q$. By Lemma~2, ${A'\leftrightarrow B'}$ is not a tautology,
which contradicts the assumption that ${A\leftrightarrow B}$ is a
tautology.

So $p$ and $q$ are directly conjunctively joined in $B$. So in
${\cal S}_{\wedge,\vee}$ we have as a theorem ${B\leftrightarrow
D}$ for a diversified formula $D$ of ${\cal L}_{\wedge,\vee}$ with
a subformula ${p\wedge q}$. Since ${A\leftrightarrow B}$ is a
tautology, it is clear that ${A\leftrightarrow D}$ is a tautology
too. If $A^{p\wedge q}_q$ and $D^{p\wedge q}_q$ are obtained from
respectively $A$ and $D$ by replacing the single occurrence of
$p\wedge q$ by $q$ , then ${A^{p\wedge q}_q\leftrightarrow
D^{p\wedge q}_q}$ is a tautology, because
${A^p_\top\leftrightarrow D^p_\top}$ is a tautology. In
$A^{p\wedge q}_q$ we have $k$ connectives, and so by the induction
hypothesis we obtain that ${A^{p\wedge q}_q\leftrightarrow
D^{p\wedge q}_q}$ is a theorem of ${\cal S}_{\wedge,\vee}$. Hence
${(A^{p\wedge q}_q)^q_{p\wedge q}\leftrightarrow (D^{p\wedge
q}_q)^q_{p\wedge q}}$ is a theorem of ${\cal S}_{\wedge,\vee}$; in
other words, ${A\leftrightarrow D}$ is a theorem of ${\cal
S}_{\wedge,\vee}$, and since ${B\leftrightarrow D}$ is such a
theorem too, we obtain that ${A\leftrightarrow B}$ is a theorem of
${\cal S}_{\wedge,\vee}$. We proceed analogously when $p\vee q$ is
a subformula of $A$.\qed

\vspace{2ex}

A $\neg\!\wedge\!\vee$\emph{-equivalence} is defined as a
$\wedge\vee$-equivalence with ${\cal L}_{\wedge\vee}$ replaced by
${\cal L}_{\neg,\wedge,\vee}$, and, as before, a
$\neg\!\wedge\!\vee$-equivalence ${A\leftrightarrow B}$ is
\emph{diversified} when $A$ and $B$ are diversified. The formal
system ${\cal S}_{\neg,\wedge,\vee}$ is defined as ${\cal
S}_{\wedge,\vee}$ save that the axioms and theorems are
$\neg\!\wedge\!\vee$-equivalences instead of
$\wedge\vee$-equivalences, and we have the additional axioms of
the following forms:
\begin{tabbing}
\centerline{${\neg\neg A\leftrightarrow A}$,}
\\[1ex]
\centerline{$\neg(A\wedge B)\leftrightarrow(\neg A\vee\neg
B)$,\hspace{6em}$\neg(A\vee B)\leftrightarrow(\neg A\wedge\neg
B)$,}
\end{tabbing}
and the following additional rule:
\[
\frac{A\leftrightarrow B}{\neg A\leftrightarrow\neg B}
\]

A formula of $\cal L$ is called $\neg$-\emph{reduced} when $\neg$
occurs in it only before letters (i.e.\ only in subformulae of the
form $\neg p$). A letter occurs \emph{positively} in a
$\neg$-reduced formula when it is not in the scope of $\neg$;
otherwise it occurs \emph{negatively}. We can prove the following.

\prop{Proposition 1$\neg$}{A diversified
$\neg\!\wedge\!\vee$-equivalence is a tautology iff it is a
theorem of ${\cal S}_{\neg,\wedge,\vee}$.}

\dkz It is clear that every theorem of ${\cal
S}_{\neg,\wedge,\vee}$ is a tautology. For the converse, we
suppose that the diversified $\neg\!\wedge\!\vee$-equivalence
${A\leftrightarrow B}$ is a tautology. It is easy to see that in
${\cal S}_{\neg,\wedge,\vee}$ we have as theorems
${A\leftrightarrow A'}$ and ${B\leftrightarrow B'}$ for
$\neg$-reduced formulae $A'$ and $B'$ of ${\cal
L}_{\neg,\wedge,\vee}$ (for that, the new axioms of ${\cal
S}_{\neg,\wedge,\vee}$ are essential), such that the diversified
$\neg\!\wedge\!\vee$-equivalence ${A'\leftrightarrow B'}$ is a
tautology.

It is impossible that a letter $p$ occurs positively in $A'$ and
negatively in $B'$. Suppose it does. Very much as in the proof of
Lemma~1, we would obtain for some $S_1,\ldots,S_n$, where
$S_i\in\{\top,\bot\}$, that $(A')^{q_1\ldots q_n}_{S_1\ldots
S_n}\leftrightarrow p$ is a tautology. On the other hand,
$(B')^{q_1\ldots q_n}_{S_1\ldots S_n}\leftrightarrow C$ can be a
tautology only if either ${C\leftrightarrow\neg p}$, or
${C\leftrightarrow\top}$, or ${C\leftrightarrow\bot}$, is a
tautology. So $(A')^{q_1\ldots q_n}_{S_1\ldots S_n}\leftrightarrow
(B')^{q_1\ldots q_n}_{S_1\ldots S_n}$ is not a tautology, and
hence ${A'\leftrightarrow B'}$ is also not a tautology,
contradicting what we have inferred above from the assumption that
${A\leftrightarrow B}$ is a tautology. It is thereby impossible
too that a letter occurs negatively in $A'$ and positively in
$B'$.

So the diversified $\neg\!\wedge\!\vee$-equivalence
${A'\leftrightarrow B'}$ is an instance of a diversified
$\wedge\vee$-equivalence ${A''\leftrightarrow B''}$, which is a
tautology. (Just replace every letter $p$ occurring negatively in
$A'$ and $B'$ by $\neg p$.) By Proposition~1, we have that
${A''\leftrightarrow B''}$ is a theorem of ${\cal
S}_{\wedge,\vee}$. Hence ${A'\leftrightarrow B'}$ is a theorem of
${\cal S}_{\neg,\wedge,\vee}$, and since ${A\leftrightarrow A'}$
and ${B\leftrightarrow B'}$ are theorems of ${\cal
S}_{\neg,\wedge,\vee}$, we obtain that ${A\leftrightarrow B}$ is a
theorem of ${\cal S}_{\neg,\wedge,\vee}$. \qed

\vspace{2ex}

It is easy to see that the formal systems ${\cal S}_{\wedge,\vee}$
and ${\cal S}_{\neg,\wedge,\vee}$ are decidable formal systems.
(For every formula of ${\cal L}_{\neg,\wedge,\vee}$ we have a
$\neg$-reduced normal form in ${\cal L}_{\neg,\wedge,\vee}$ unique
up to associativity and commutativity of $\wedge$ and $\vee$.)

\section{Isomorphic formulae with perfect generalizability}
Let $\cal K$ be a category whose objects are the formulae of
${\cal L}_{\neg,\wedge,\vee}$, and whose arrows $f\!:A\rightarrow
B$ are intuitively interpreted as deductions, or proofs, from $A$
to $B$. We assume that $\cal K$ has isomorphisms covering the
theorems of ${\cal S}_{\neg,\wedge,\vee}$. This means, for
example, that we have an isomorphism of the type $(A\wedge
B)\wedge C\rightarrow A\wedge(B\wedge C)$, whose inverse is of the
type $A\wedge(B\wedge C)\rightarrow(A\wedge B)\wedge C$.

These isomorphisms of $\cal K$ correspond to deductions in the
multiplicative fragment of linear propositional logic. These are
the basic arrows of $\cal K$, but $\cal K$ can have more arrows
than that. We assume however that $\cal K$ does not go beyond
deductions of classical propositional logic, which means that if
$f\!:A\rightarrow B$ is an arrow of $\cal K$, then $A\rightarrow
B$ is a tautology, with the arrow $\rightarrow$ interpreted as the
connective of material implication.

For a formula $A$, let $|A|$ be the number of occurrences of
letters in $A$. For every arrow $f\!:A\rightarrow B$ of $\cal K$,
if $x$ is an occurrence of a letter in $A$, then let $o(x)=n\mn 1$
iff $x$ is the $n$-th occurrence of letter in $A$ counting from
the left, and if $x$ is an occurrence of a letter in $B$, then let
$o(x)=|A|\pl n\mn 1$ iff $x$ is the $n$-th occurrence of letter in
$B$ counting from the left.

We assume that every arrow $f\!:A\rightarrow B$ of $\cal K$
induces on the ordinal $|A| + |B|$ an equivalence relation $L_f$,
called the \emph{linking relation} of $f$, which satisfies the
condition that $(o(x),o(y))\in L_f$ \emph{only if} $x$ and $y$ are
occurrences in $A$ or $B$ of the same letter. A linking relation
is \emph{perfect} when instead of \emph{only if} in this condition
we have \emph{if and only if}.

We say that $A$ and $B$ are \emph{uniform instances} of $A_1$ and
$B_1$ when they are instances of $A_1$ and $B_1$ respectively with
the same letter-for-letter substitution (i.e.\ substitution that
replaces a letter by a letter).

We say that an arrow $f\!:A\rightarrow B$ is \emph{generalized} to
an arrow $f_1\!:A_1\rightarrow B_1$ when $A$ and $B$ are uniform
instances of $A_1$ and $B_1$, and the linking relations $L_f$ and
$L_{f_1}$ are the same.

We say that a category $\cal K$ is \emph{perfectly generalizable}
when each of its arrows can be generalized to an arrow of $\cal K$
with a perfect linking relation. Examples of perfectly
generalizable categories with syntactically defined arrows will be
given below after Lemma~3.

A category $\cal K$ will be called \emph{permutational} when for
every isomorphism $f\!:A\rightarrow B$ of $\cal K$ the linking
relation $L_f$ corresponds to a bijection between $|A|$ and $|B|$,
and if $g\!:B\rightarrow A$ is the inverse of $f$, then $L_g$
corresponds to the inverse of this bijection. Permutational
categories may be both perfectly generalizable and not perfectly
generalizable. In this section we consider those of the first
kind, while in the next one we deal with the second kind.

Permutational categories arise naturally whenever we have a
certain modelling of categories that we are now going to describe.
The category \emph{SplPre} is the category whose arrows are
\emph{split preorders} between finite ordinals; \emph{Gen} is a
subcategory of \emph{SplPre} whose arrows are \emph{split
equivalences} between finite ordinals, while \emph{Rel} is a
category whose arrows are relations between finite ordinals. (The
category \emph{Rel} has an isomorphic image within \emph{SplPre}
by a map that preserves composition, but not identity arrows.) A
split preorder is a preordering relation on the disjoint union of
two sets conceived as source and target, and analogously for split
equivalences, which are equivalence relations. We have
investigated \emph{SplPre} and \emph{Gen} systematically in
\cite{DP09}, \cite{DP03a} and \cite{DP03b}, and we have used
\emph{Gen} and \emph{Rel} as model categories for equality of
deductions in \cite{DP04} and \cite{DP07}, and in other papers
related to these two books.

We may use these model categories to produce the linking relation
of $\cal K$ in the following manner. For a functor $G$ from $\cal
K$ into a model category such as \emph{SplPre}, \emph{Gen} or
\emph{Rel}, we take that $(n,m)\in L_f$ iff the ordinals
corresponding to $n$ and $m$ are linked by the reflexive,
symmetric and transitive closure of $Gf$, which coincides with
$Gf$ if $Gf$ is an equivalence relation of \emph{Gen}. In this
section we are interested in particular in linking relations
produced by \emph{Gen}, while in the next section we will
encounter also one produced by \emph{Rel}.

It may happen that the same category $\cal K$ produces different
linking relations with different functors $G$. The notions of
perfectly generalizable and permutational category are relative to
a chosen kind of linking relations.

Let us give an example of linking relations. We may have in $\cal
K$ an arrow $f\!:p\wedge(\neg p\vee p)\rightarrow p$ corresponding
to modus ponens, such that $L_f$ will give the following linking

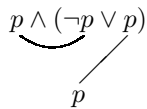
\begin{center}
\begin{picture}(50,50)

\qbezier(3,34)(15,25)(27,34)

\put(25.5,16){\line(1,1){18}}


\put(25,10){\makebox(0,0){$p$}}

\put(25,40){\makebox(0,0){$p\wedge (\neg p\vee p)$}}

\end{picture}
\end{center}
and another arrow $g\!:p\wedge(\neg p\vee p)\rightarrow p$
corresponding to the first projection, such that $L_f$ will give
the following linking

\begin{center}
\begin{picture}(50,50)

\put(23,16){\line(-1,1){18}}

\put(25,10){\makebox(0,0){$p$}}

\put(25,40){\makebox(0,0){$p\wedge (\neg p\vee p)$}}

\end{picture}
\end{center}
The arrow $f$ can be generalized to an arrow of the type
$p\wedge(\neg p\vee q)\rightarrow q$, while $g$ can be generalized
to an arrow of the type $p\wedge(\neg q\vee r)\rightarrow p$. Both
of these arrows to which $f$ and $g$ are generalized have a
perfect linking.

If the linking relations of the arrows of $\cal K$ are produced by
a functor $G$ into \emph{SplPre} (or a subcategory thereof), as
explained above, then we may conclude that $\cal K$ is
permutational. This is because the isomorphisms of \emph{SplPre}
correspond to bijections. This is so both for perfectly
generalizable and for not perfectly generalizable categories $\cal
K$. We can prove the following lemma.

\prop{Lemma~3}{If $\cal K$ is a permutational perfectly
generalizable category, and $A$ and $B$ are isomorphic in $\cal
K$, then there are diversified formulae $A_1$ and $B_1$ such that
$A$ and $B$ are uniform instances of $A_1$ and $B_1$ and
${A_1\leftrightarrow B_1}$ is a tautology.}

\dkz Since $\cal K$ is perfectly generalizable, the arrows
$f\!:A\rightarrow B$ and its inverse $g\!:B\rightarrow A$ of $\cal
K$ are generalized to the arrows $f_1\!:A_1\rightarrow B_1$ and
$g_2\!:B_2\rightarrow A_2$ with perfect linking relations. Since
$\cal K$ is permutational, these relations correspond to
bijections inverse to each other. From that we conclude that
$A_1$, $B_1$, $A_2$ and $B_2$ are diversified formulae. Since $B$
is a letter-for-letter instance of both $B_1$ and $B_2$, which are
diversified, $B_1$ is a letter-for-letter instance of $B_2$, and
since the linking relations of $f_1$ and $g_2$ correspond to
bijections inverse to each other, our letter-for-letter
substitution produces $A_1$ out of $A_2$. According to what we
assumed at the beginning of the section, $A_1\rightarrow B_1$ and
$B_2\rightarrow A_2$ are tautologies, and hence $B_1\rightarrow
A_1$ is a tautology too. \qed

\vspace{2ex}

As categories $\cal K$ covered by this lemma we have the free
distributive lattice category of \cite{DP04} (Section 11.1), which
axiomatizes equality of deductions in conjunctive-disjunctive
classical logic (the objects of this category are in ${\cal
L}_{\wedge,\vee}$), and the free proof-net category of \cite{DP07}
(Section 2.2), which axiomatizes equality of deductions in the
multiplicative fragment of linear logic without propositional
constants (its objects are in ${\cal L}_{\neg,\wedge,\vee}$). One
may easily conceive other such examples, and in particular the
example of a category axiomatizing equality of deductions in the
whole of classical propositional logic (whose objects are in
${\cal L}_{\neg,\wedge,\vee}$). One would just take as equations
of $\cal K$ those equations $f=g$ such that $Gf=Gg$, where $G$ is
a nontrivial functor from $\cal K$ to \emph{Gen} making $\cal K$
perfectly generalizable.

The following proposition characterizes pairs of isomorphic
formulae for permutational perfectly generalizable categories.

\prop{Proposition 2$\neg$}{The formulae $A$ and $B$ are isomorphic
in a permutational perfectly generalizable category iff
${A\leftrightarrow B}$ is a theorem of ${\cal
S}_{\neg,\wedge,\vee}$.}

\dkz From left to right, suppose that $A$ is isomorphic to $B$.
So, by Lemma~3, there are diversified formulae $A_1$ and $B_1$
such that $A$ and $B$ are uniform instances of $A_1$ and $B_1$,
and ${A_1\leftrightarrow B_1}$ is a tautology. By
Proposition~1$\neg$, we obtain that ${A_1\leftrightarrow B_1}$ is
a theorem of ${\cal S}_{\neg,\wedge,\vee}$, and hence
${A\leftrightarrow B}$ is that too.

For the other direction, we have assumed that the equivalences of
${\cal S}_{\neg,\wedge,\vee}$ are covered by the isomorphisms of
our category. \qed

\vspace{2ex}

\noindent We have a Proposition~2 analogous to Proposition~2$\neg$
with ${\cal S}_{\neg,\wedge,\vee}$ replaced by ${\cal
S}_{\wedge,\vee}$, for categories whose objects are the formulae
of ${\cal L}_{\wedge,\vee}$. Since both ${\cal S}_{\wedge,\vee}$
and ${\cal S}_{\neg,\wedge,\vee}$ are easily seen to be decidable
systems, Propositions 2 and 2$\neg$ lead to a decision procedure
for isomorphic formulae.

\section{Isomorphic formulae without perfect generalizability}
For $f\!:A\rightarrow B$ an arrow, let the set of \emph{diagonal
links} of $f$ be the set $D_f=\{(n,n)\mid n\in |A|\pl|B|\}$. We
will here call an arrow $f\!:A\rightarrow B$ in a category a
\emph{zero arrow} when its linking relation $L_f$ is equal to the
set $D_f$ of its diagonal links. A \emph{zero identity} arrow is
an arrow $0_A\!:A\rightarrow A$ with ${L_{0_A}\!=D_{0_A}}$. The
\emph{union} of the arrows $f,g\!:A\rightarrow B$ will here be an
arrow $f\cup g\!:A\rightarrow B$ such that $L_{f\cup g}$ is the
transitive closure of $L_f\cup L_g$. (For a general treatment of
zero arrows and union of arrows see \cite{DP04}.)

Consider a category $\cal K$ such that the objects of $\cal K$ are
the formulae of $\cal L$ and its arrows correspond to deductions
in classical propositional logic. We assume that in $\cal K$ we
have for every $A$ of $\cal L$ a zero identity arrow $0_A$, and we
have also closure under union of arrows. As a particular case of
$\cal K$ we have the category \textbf{B} of \cite{DP04} (Section
14.2), but possibly also a category that besides the zero arrows
of the types $\top\rightarrow \neg A\vee A$ and $A\wedge\neg
A\rightarrow\bot$, like those we have in \textbf{B}, has also
arrows of this type with linking relations given by the following
pictures:

\begin{center}
\begin{picture}(150,43)

\qbezier(17,16)(26,25)(35,16)

\put(22,10){\makebox(0,0){$\neg A\vee A$}}

\put(25,33){\makebox(0,0){$\top$}}

\qbezier(115,28)(127,19)(140,28)

\put(128,33){\makebox(0,0){$A\wedge\neg A$}}

\put(125,10){\makebox(0,0){$\bot$}}

\end{picture}
\end{center}

The linking relation $L_f$ of an arrow $f\!:A\rightarrow B$ of
\textbf{B}, which is an equivalence relation, is not the relation
$Gf\subseteq {A\times B}$, where $G$ is the functor into
\emph{Rel} with respect to which \textbf{B} is coherent; i.e., the
functor $G$ such that for $f,g\!:A\rightarrow B$ arrows of
\textbf{B} we have that $f=g$ iff $Gf=Gg$ (see \cite{DP04},
Section 14.2). We define $L_f$ as the reflexive, symmetric and
transitive closure of $Gf$. Note that for $f,g\!:A\rightarrow B$
arrows of \textbf{B} we do not have that $L_f=L_g$ only if $f=g$
in \textbf{B}, but we still have that for particular arrows $f$
and $g$, as in the proof of Lemma~7. What we need essentially is
that this holds when one of $f$ and $g$ is an identity arrow.

For the categories $\cal K$ other than \textbf{B} we need the same
assumption, which is guaranteed of course if equality of arrows in
$\cal K$ is defined via the linking relations---if we have namely
that for $f,g\!:A\rightarrow B$ the equation $f=g$ holds in $\cal
K$ iff $L_f=L_g$. Our assumptions about $\cal K$ must guarantee
also that we have an analogue of Lemma~6 below. These assumptions
can be quite standard, like those spelled out in \cite{DP04} and
\cite{DP07}. If composition of these linking relations is defined
naturally, like composition in the category \emph{Gen} (see
\cite{DP09}, Section~2, for a detailed study), then in \textbf{B}
we do not have that $L_{g\cirkk f}$ is equal to the composition of
$L_f$ with $L_g$, but this may hold for other categories $\cal K$.

We take that our assumptions imply that $\cal K$ is permutational
in the sense of the preceding section. However, $\cal K$ is not
perfectly generalizable, because of the presence of zero identity
arrows. The arrow $0_p\!:p\rightarrow p$ cannot be generalized to
one with a perfect linking relation. The category \textbf{B} is
permutational and is not perfectly generalizable.

The definition of the linking relation for the arrows of
\textbf{B} given above, and our notion of linking relation for
$\cal K$ in general, is motivated by the wish to have a uniform
notion based, as in the preceding section, on an equivalence
relation. This uniform notion is related to generality of
deductions (more consistently in the preceding section than in the
present one). If however we do not strive for this uniformity,
then we may define $L_f$ for \textbf{B} just as $Gf$, and our
proofs would still go through. For $\cal K$ in general, we may
also have a notion of linking relation like $Gf$.

We prove first a simple preliminary lemma.

\prop{Lemma~4}{For every formula $A$ of ${\cal
L}_{\top,\bot,\wedge,\vee}$ in which the letter $p$ occurs there
are formulae $A_1$ and $A_2$ of ${\cal L}_{\top,\bot,\wedge,\vee}$
such that $A\leftrightarrow((p\wedge A_1)\vee A_2)$ is a
tautology.}

\dkz We proceed by induction on the number $n$ of occurrences of
binary connectives in $A$. If $n=0$, then $A$ is $p$, and we take
$A_1$ and $A_2$ to be respectively $\top$ and $\bot$.

Suppose $n>0$. If $A$ is $A'\wedge A''$ with an occurrence of $p$
in $A'$, then by the induction hypothesis we have formulae $A_1'$
and $A_2'$ of ${\cal L}_{\top,\bot,\wedge,\vee}$ such that
$A'\leftrightarrow((p\wedge A_1')\vee A_2')$ is a tautology. We
take $A_1$ and $A_2$ to be respectively $A_1'\wedge A''$ and
$A_2'\wedge A''$. If $A$ is $A'\vee A''$ with an occurrence of $p$
in $A'$, then by the induction hypothesis we have formulae $A_1'$
and $A_2'$ as before, and we take $A_1$ and $A_2$ to be
respectively $A_1'$ and $A_2'\vee A''$. If $A$ is $A''\wedge A'$
or $A''\vee A'$ with an occurrence of $p$ in $A'$, then we have
that $A\leftrightarrow(A'\wedge A'')$ or $A\leftrightarrow(A'\vee
A'')$ is a tautology, and we proceed as before. \qed

\vspace{2ex}

\noindent The tautology of this lemma is provable in every logic
that algebraically corresponds to distributive lattices with top
and bottom. We have next the following lemmata.

\prop{Lemma~5}{If $B'$ is obtained from a formula $B$ of ${\cal
L}_{\top,\bot,\wedge,\vee}$ by replacing a single occurrence of
$q$ by $p\wedge q$, then $(p\wedge B)\rightarrow B'$ is a
tautology.}

\dkz We proceed by induction on he number $n$ of occurrences of
binary connectives in $B$. If $n=0$, then $B$ is $q$, and
$(p\wedge q)\rightarrow(p\wedge q)$ is a tautology.

Suppose $n>0$. If $B$ is $B_1\wedge B_2$ with the replaced
occurrence of $q$ in $B_1$, then by the induction hypothesis we
have that $(p\wedge B_1)\rightarrow B_1'$ is a tautology, and
hence, by relying on the associativity of $\wedge$, we have that
$(p\wedge(B_1\wedge B_2))\rightarrow B'$, where $B'$ is
$B_1'\wedge B_2$, is a tautology too. If $B$ is $B_1\vee B_2$ with
the replaced occurrence of $q$ in $B_1$, by the induction
hypothesis we have again that $(p\wedge B_1)\rightarrow B_1'$ is a
tautology, and since $(p\wedge(B_1\vee B_2))\rightarrow((p\wedge
B_1)\vee B_2)$ (which corresponds to a dissociativity arrow; see
\cite{DP04}, Chapter~7) is a tautology, we have that
$(p\wedge(B_1\vee B_2))\rightarrow B'$, where $B'$ is $B_1'\vee
B_2$, is a tautology too. If $B$ is $B_2\wedge B_1$ or $B_2\vee
B_1$ with the replaced occurrence of $q$ in $B_1$, then we have
that $B\rightarrow(B_1\wedge B_2)$ or $B\rightarrow(B_1\vee B_2)$
is a tautology, and we proceed as before. \qed

\prop{Lemma~6}{If $A$ and $B$ are formulae of ${\cal
L}_{\top,\bot,\wedge,\vee}$, with $x$ and $y$ occurrences of the
letter $p$ in respectively $A$ and $B$, and $A\rightarrow B$ is a
tautology, then there is an arrow $f\!:A\rightarrow B$ of
\mbox{\rm \textbf{B}} such that in ${L_f-D_f}$ we find only the
pairs $(o(x),o(y))$ and $(o(y),o(x))$.}

\dkz By Lemma~4, we have formulae $A_1$ and $A_2$ of ${\cal
L}_{\top,\bot,\wedge,\vee}$ such that $A\leftrightarrow((p\wedge
A_1)\vee A_2)$ is a tautology. The proof of this lemma yields an
arrow ${\tau\!:A\rightarrow(p\wedge A_1)\vee A_2}$ such that
$(o(x),|A|)\in L_\tau$ (note that if $z$ is the leftmost
occurrence of $p$ in $(p\wedge A_1)\vee A_2$, i.e.\ the occurrence
written down, then $o(z)=|A|$). We also have arrows
$\sigma\!:(p\wedge A_1)\vee A_2\rightarrow A$ and
$g\!:A\rightarrow B$ of \textbf{B}, because their types correspond
to tautologies. If $\iota_1\!:p\wedge A_1\rightarrow (p\wedge
A_1)\vee A_2$ and $\iota_2\!:A_2\rightarrow(p\wedge A_1)\vee A_2$
are injection arrows, then for
\[
\zeta_i=_{df} 0_B\cirk g\cirk\sigma\cirk\iota_i,
\]
where $i\in\{1,2\}$, we have $L_{\zeta_i}=D_{\zeta_i}$, because
the composition ends with $0_B$. If $\pi\!:p\wedge A_1\rightarrow
p$ is a first projection arrow, then for
$\langle\pi,\zeta_1\rangle\!:p\wedge A_1\rightarrow p\wedge B$ we
have
\[
L_{\langle\pi,\zeta_1\rangle}-D_{\langle\pi,\zeta_1\rangle}=\{(0,|A_1|\pl
1),(|A_1|\pl 1,0)\}
\]
(this means that the two occurrences of $p$ that are written down
are linked).

If $B'$ is obtained from $B$ by replacing the occurrence $y$ of
$p$ in $B$ by $p\wedge p$, then the proof of Lemma~5 yields an
arrow $\eta\!:p\wedge B\rightarrow B'$ of \textbf{B}. Out of the
first projection arrow $\pi'\!:p\wedge p\rightarrow p$ we obtain
an arrow $\theta\!:B'\rightarrow B$ such that for
$\theta\cirkk\eta\!:p\wedge B\rightarrow B$ we have $(0,o(y))\in
L_{\theta\cirk\eta}$ (i.e.\ the occurrence of $p$ written down in
the source $p\wedge B$ and the occurrence $y$ in the target $B$
are linked).

For
\[
\mu=_{df}[\theta\cirk\eta\cirk\langle\pi,\zeta_1\rangle,\zeta_2]\!:(p\wedge
A_1)\vee A_2\rightarrow B
\]
we have $L_\mu-D_\mu=\{(0,o(y)),(o(y),0)\}$, and we take
$f\!:A\rightarrow B$ to be $\mu\cirk\tau$. We have
$L_{\mu\cirkk\tau}-D_{\mu\cirkk\tau}=\{(o(x),o(y)),(o(y),o(x))\}$.
\qed

\vspace{2ex}

The formulae $A$ and $B$ of ${\cal L}_{\top,\bot,\wedge,\vee}$ are
\emph{letter-homogeneous} when every letter that occurs in $A$
occurs in $B$ the same number of times, and vice versa, every
letter that occurs in $B$ occurs in $A$ the same number of times.
Two $\neg$-reduced formulae of $\cal L$ are letter-homogeneous
when they are uniform instances of two letter-homogeneous formulae
of ${\cal L}_{\top,\bot,\wedge,\vee}$. This means that such
formulae share both positive and negative occurrences of letters.
We can prove the following.

\prop{Lemma~7}{For $A$ and $B$ formulae of ${\cal
L}_{\top,\bot,\wedge,\vee}$, if ${A\leftrightarrow B}$ is a
tautology and $A$ and $B$ are letter-homogeneous, then $A$ and $B$
are isomorphic in \mbox{\rm \textbf{B}}.}

\dkz Suppose $A\rightarrow B$ is a tautology and $A$ and $B$ are
letter-homogeneous. For every bijection mapping the occurrences of
letters in $A$ to the occurrences of the same letters in $B$,
there is an arrow $f\!:A\rightarrow B$ of \textbf{B} such that
$L_f$ corresponds to this bijection. This is guaranteed by Lemma~6
and the operation of union of arrows of \textbf{B}. Since the same
holds for the tautology $B\rightarrow A$ and the inverse of our
bijection, we conclude that there is an arrow $g\!:B\rightarrow A$
of \textbf{B} such that $L_g$ corresponds to this inverse. With
the help of Boolean Coherence (see \cite{DP04}, Section 14.2; this
is the assertion that the functor $G$, on which the definition of
$L_f$ is based, is a faithful functor from \textbf{B} to
\emph{Rel}), we infer from $L_{g\cirkk f}=L_{\mj_A}$ and
$L_{f\cirkk g}=L_{\mj_B}$ that $f$ and $g$ are isomorphisms,
inverse to each other. \qed

\vspace{2ex}

We can now prove the following.

\prop{Proposition~3}{For $A$ and $B$ $\neg$-reduced formulae of
$\cal L$ we have that $A$ and $B$ are isomorphic in \mbox{\rm
\textbf{B}} iff ${A\leftrightarrow B}$ is a tautology and $A$ and
$B$ are letter-homogeneous.}

\dkz The direction from left to right is an easy consequence of
the fact that the arrows of \textbf{B} correspond to implications
that are tautologies, of the fact that no arrow of \textbf{B}
links a positive occurrence of a letter in the source with a
negative occurrence of a letter in the target, and vice versa, and
of the fact that \textbf{B} is a permutational category in the
sense of Section~3.

For the other direction we use Lemma~7 and appeal to the fact that
if the right-hand side holds, then there are formulae $A_1$ and
$B_1$ of ${\cal L}_{\top,\bot,\wedge,\vee}$ such that $A_1$ and
$B_1$ are letter-homogeneous, $A$ and $B$ are uniform instances of
$A_1$ and $B_1$ and ${A_1\leftrightarrow B_1}$ is a tautology. We
derive that ${A_1\leftrightarrow B_1}$ is a tautology by
appropriate substitutions within ${A\leftrightarrow B}$. \qed

\vspace{2ex}

Every formula of $\cal L$ may effectively be reduced to a
$\neg$-reduced formula isomorphic in \textbf{B}. So Proposition~3
gives a characterization of arbitrary pairs of isomorphic formulae
of \textbf{B}, and leads to a decision procedure for checking
whether such an isomorphism exists.

\vspace{4ex}

\noindent {\small {\it Acknowledgement.} Work on this paper was
supported by the Ministry of Science of Serbia (Grants 144013 and
144029).}

\end{document}